\documentclass{article}
\usepackage{amssymb}
\usepackage{amsmath}
\usepackage{amsthm}
\usepackage{bbold}

\newcounter{enucount}
\newcounter{enuref}
\newcounter{enualph}
\newcounter{enunumb}
\renewcommand{\theenucount}{(\roman{enucount})}

\renewcommand{\theenunumb}{(\arabic{enunumb})}
\renewcommand{\theenualph}{(\alph{enualph})}

\theoremstyle{plain}
\newtheorem{theorem}{Theorem}

\newtheorem{lem}[theorem]{Lemma}

\newtheorem{fact}[theorem]{Fact}

\newtheorem{sclaim}{Claim}[theorem]

\theoremstyle{definition}
\newtheorem{defin}{Definition}
\newtheorem{exam}{Example}

\theoremstyle{remark}

\renewenvironment{enumerate}{\begin{list}{\rm \theenucount}{\usecounter{enucount}
    \setlength{\labelwidth}{1cm}}}
   {\end{list}}

\newcommand{\embed}{<\!\!\circ\;}

\newcommand{\D}{{\cal D}}

\renewcommand{\P}{{\cal P}}

\newcommand{\U}{{\cal U}}

\renewcommand{\AA}{{\mathbb A}}
\newcommand{\BB}{{\mathbb B}}
\newcommand{\CC}{{\mathbb C}}

\newcommand{\EE}{{\mathbb E}}

\newcommand{\PP}{{\mathbb P}}
\newcommand{\QQ}{{\mathbb Q}}

\renewcommand{\aa}{{\mathfrak a}}
\newcommand{\bb}{{\mathfrak b}}
\newcommand{\cc}{{\mathfrak c}}
\newcommand{\dd}{{\mathfrak d}}

\newcommand{\sss}{{\mathfrak s}}

\newcommand{\uu}{{\mathfrak u}}

\newcommand{\cov}{{\mathsf{cov}}}
\newcommand{\non}{{\mathsf{non}}}

\newcommand{\amal}{{\mathrm{amal}}}

\newcommand{\ro}{{\mathrm{r.o.}}}
\newcommand{\dir}{{\mathrm{dir}}}
\newcommand{\cl}{{\mathrm{cl}}}

\newcommand{\zero}{{\mathbb{0}}}
\newcommand{\one}{{\mathbb{1}}}

\newcommand{\sub}{\subseteq}
\newcommand{\sem}{\setminus}

\newcommand{\la}{\langle}
\newcommand{\ra}{\rangle}

\newcommand{\noi}{\noindent}

\title{The amalgamated limit and its topological interpretation}

\author{J\"org Brendle\thanks{Partially supported by Grant-in-Aid for Scientific Research
   (C) 18K03398, Japan Society for the Promotion of Science. This survey is related to the author's minicourse {\em Recent forcing
   techniques: Boolean ultrapowers} at the RIMS workshop {\em Recent Developments in Set Theory of the Reals}
   in October 2021. He thanks the Research Institute for Mathematical Sciences (RIMS), Kyoto University, for their
   support and Masaru Kada for inviting him. \newline
       }  \\ 
   Graduate School of System Informatics \\
   Kobe University \\
   Rokko-dai 1-1, Nada-ku \\
   Kobe 657-8501, Japan \\   
   email: {\sf brendle@kobe-u.ac.jp}}

\begin{document}
\maketitle

\begin{abstract}
\noi This is a survey on the amalgamated limit from~\cite{shattered}, a limit construction for complete Boolean algebras in iterated forcing theory,
which generalizes both the direct limit and the two-step amalgamation. We focus in particular on examples of the amalgamated limit
from the literature and on the topological amalgamated limit for compact Hausdorff spaces.
\end{abstract}



\section*{Introduction}

Limit constructions for complete Boolean algebras (cBa's for short) play a fundamental role in forcing theory because they are needed to
set up iterated forcing constructions. The most important limits certainly are the direct limit and the inverse limit and both have been used for decades
in a myriad of consistency proofs. The present note deals with a construction generalizing both the direct limit and the two-step 
amalgamation, the {\em amalgamated limit}, which has been introduced in~\cite{shattered} for ccc iterations adding random reals instead
of Cohen reals in limit stages, called {\em shattered iterations} (see Example~\ref{exam5} in Section 2). However, at least implicitly, the amalgamated limit has been used
much earlier in work of Shelah (\cite{Sh700}, see also~\cite{Br07}, \cite[Section 1]{survey}, and Example~\ref{exam4} in Section 2). Roughly
speaking, the amalgamated limit $\AA_\ell = \lim\amal_{i \in I} \AA_i$ of a system of cBa's $(\AA_i : i \in I)$ is the natural ``smallest" cBa
into which all $\AA_i$ completely embed where $I$ is a distributive almost-lattice (see Section 1 for a formal definition); the latter can be a fairly general
structure, but for applications one can think of $I$ as being a set of pairs of ordinals. Note the similarity to the direct limit here: the latter is
the smallest cBa completely embedding all $\AA_i$ when $I$ is a lattice -- or just a directed set. Also notice that since distributive almost-lattices
are inherently two-dimensional objects, there is a close connection of the amalgamated limit with {\em matrix iterations} which will be made precise
in Example~\ref{exam3} in Section 2. However, since matrix iterations work with finite support (i.e., direct limits are taken everywhere),
the amalgamated limit framework is not necessary for them and in a sense just complicates things.

The present note is a survey on the amalgamated limit. Most of the material presented here has appeared (or will appear) elsewhere,
mostly in~\cite{shattered}, and to a lesser extent in~\cite{survey}. For this reason we refrain from giving proofs and just present the main notions
and ideas. Section 1 contains the basics: correctness (originally from~\cite{Br05}, see also~\cite[Section 1]{shattered} and~\cite{survey}) is reviewed,
the amalgamated limit is defined, and its most important properties (in particular, complete embeddability, Theorem~\ref{limamal-cBa}) are stated.
We kept this section brief because all the material can be found in~\cite{shattered}. Section 2 contains the most important examples we are aware
of, and which are scattered through the literature. Again, we do not present concrete and detailed proofs (because they can be found elsewhere)
but rather give the main ideas of the constructions and explain why they can be construed as  amalgamated limits. Section 3, finally, is the
only section with completely new material: by Stone duality, the amalgamated limit can be redone for duals of cBa's  -- or actually, Boolean algebras
in general -- to yield a limit construction for compact zero-dimensional Hausdorff spaces. It turns out that zero-dimensionality is irrelevant here and
that this construction even works for general compact Hausdorff spaces. The limit space $X_\ell$ is a subspace of the product space of the $X_i$, $
i \in I$, where $I$ is again a distributive almost-lattice, such that the projection mappings are open (see Theorem~\ref{limamal-top}). 
The latter -- corresponding to complete embeddability in the cBa setting of Theorem~\ref{limamal-cBa} -- is the main point. In this last section we present detailed proofs.

We assume familiarity with forcing theory, in particular with iterated forcing constructions and with the approach to forcing via cBa's.
Basics can be found in the standard literature, e.g.~\cite{Je03}, \cite{Ku11}, and~\cite{Ha17}. Furthermore we will mention a number
of standard cardinal invariants of the continuum along the way without definition, see~\cite{BJ95}, \cite{Bl10}, or~\cite{Ha17}.


\section{The amalgamated limit}

For complete Boolean algebras $\AA_i \sub \AA_j$, we write $\AA_i \embed \AA_j$ for ``{\em $\AA_i$ is completely
embedded in $\AA_j$}". Recall that this is equivalent to saying that the {\em projection mapping} $h^j_i : \AA_j \to \AA_i$ given by
\[ h^j_i (a_j) = \prod \{ a_i \in \AA_i : a_i \geq a_j \} \]
satisfies $h^j_i (a_j) > \zero$ for all $a_j > \zero$. We need the following basic definition from~\cite{Br05} (see also~\cite{survey} or~\cite[Definition 1]{shattered}).

\begin{defin} \label{correct-def}
Assume we have cBa's $\AA_{0 \land 1} \embed \AA_i \embed \AA_{0 \lor 1}$, $i \in \{ 0,1\}$. We say {\em projections in the diagram 
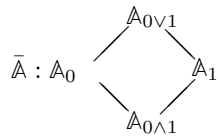
\[ 
\begin{picture}(80,50)(0,0)
\put(48,5){\makebox(0,0){$\AA_{0\land 1}$}}
\put(7,25){\makebox(0,0){$\bar \AA :  \AA_{0}$}}
\put(68,25){\makebox(0,0){$\AA_{1}$}}
\put(48,45){\makebox(0,0){$\AA_{0\lor 1}$}}
\put(26,28){\line(1,1){13}}
\put(28,20){\line(1,-1){11}}
\put(51,8){\line(1,1){13}}
\put(53,40){\line(1,-1){11}}
\end{picture}
\]
are correct} if either of the following three equivalent conditions holds:
\begin{itemize}
\item $h^{0 \lor 1}_1 (a_0) = h^0_{0 \land 1} (a_0)$ for all $a_0 \in \AA_0$,
\item $h^{0 \lor 1}_0 (a_1) = h^1_{0 \land 1} (a_1)$ for all $a_1 \in \AA_1$,
\item whenever $h^0_{0 \land 1} (a_0) = h^1_{0 \land 1} (a_1)$ then $a_0$ and $a_1$
are compatible in $\AA_{0 \lor 1}$.  \hfill $\dashv$
\end{itemize}
\end{defin}

Notice this implies (but is not equivalent to)  $\AA_{0\land 1} = \AA_0 \cap \AA_1$. 
A typical example for a diagram with correct projections is given by letting $\AA_{0\land 1} = \{ \zero, \one \}$ and
$\AA_{0 \lor 1}$ the usual product forcing, that is, the completion of $(\AA_0 \sem \{ \zero \} ) \times (\AA_1 \sem \{ \zero \})$.
Another important example is obtained by letting $\AA_{0 \land 1} \embed \AA_0$ be arbitrary forcing notions and putting
$\AA_1 := \AA_{0\land 1} \star \dot \QQ$ and $\AA_{0 \lor 1} := \AA_0 \star \dot \QQ$, where $\QQ$ is a Suslin ccc forcing notion
(see~\cite[Example 3]{shattered}).
In both cases, correctness is straightforward. For an example of a non-correct diagram with $\AA_{0\land 1} = \AA_0 \cap \AA_1$ 
see~\cite[Counterexample 4]{shattered}.  More on correctness can be found in Section 1 of the latter work.

Since dealing with limit constructions is one of the main goals of this work, we recall:

\begin{lem}[{embeddability of direct limits, see~\cite[Lemma 4]{survey} and~\cite[Lemma 4]{shattered}}]   \label{correctness-limdir}
Let $K$ be a directed index set. Assume $( \AA_k : k \in K )$ and $( \EE_k : k \in K )$ are systems of cBa's such that $\AA_k \embed \AA_\ell$,
$\EE_k \embed \EE_\ell$, and $\AA_k \embed \EE_k$ for any $k \leq \ell$. Assume further projections in all diagrams of the form
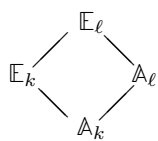
\[ 
\begin{picture}(80,50)(0,0)
\put(48,5){\makebox(0,0){$\AA_k$}}
\put(22,25){\makebox(0,0){$\EE_k$}}
\put(68,25){\makebox(0,0){$\AA_\ell$}}
\put(48,45){\makebox(0,0){$\EE_\ell$}}
\put(26,28){\line(1,1){13}}
\put(28,20){\line(1,-1){11}}
\put(51,8){\line(1,1){13}}
\put(53,40){\line(1,-1){11}}
\end{picture}
\]
are correct for $k \leq \ell$. Then $\AA : = \lim \dir_{k \in K} \AA_k$ completely embeds into $\EE := \lim \dir_{k \in K} \EE_k$. 
Furthermore, correctness is preserved in the sense that projections in all diagrams of the form
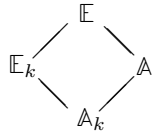
\[ 
\begin{picture}(80,50)(0,0)
\put(48,5){\makebox(0,0){$\AA_k$}}
\put(22,25){\makebox(0,0){$\EE_k$}}
\put(68,25){\makebox(0,0){$\AA$}}
\put(46,45){\makebox(0,0){$\EE$}}
\put(26,28){\line(1,1){13}}
\put(28,20){\line(1,-1){11}}
\put(51,8){\line(1,1){13}}
\put(53,40){\line(1,-1){11}}
\end{picture}
\]
for $k \in K$ are correct.
\end{lem}

Correctness is crucial here. For a (non-correct) counterexample see~\cite[Counterexample 5]{shattered}. The notion of correctness
was originally singled out by the author~\cite{Br05}  when he was grappling with Shelah's theory of iterations along templates (\cite{Sh700},
see also~\cite{Br02} and~\cite{Br03}). See~\cite{Br05} or~\cite{survey} for how template iterations can be nicely described in the
correctness framework.
 
We next define the structure underlying the amalgamated limit (see Definition~\ref{amallim-def} below).

\begin{defin}[{see~\cite[Definition 2]{shattered}}]  \label{almostlattice-def}
Assume $\la I, \leq \ra$ is a partial order. Call $I$ a {\em distributive almost-lattice} if
\begin{enumerate}
\item any two elements $i,j \in I$ have a greatest lower bound $i \land j$, the {\em meet}
   of $i$ and $j$,
\item if $i$ and $j$ have an upper bound, then they have a least upper bound $i \lor j$, the
   {\em join} of $i$ and $j$,
\item if $i_0, i_1, i_2 \in I$ then there are $j \neq k$ such that $i_j $ and $i_k$ have an upper bound
   (so that $i_j \lor i_k$ exists), that is, given any three elements, two of them have an upper bound,
\item the distributive laws hold for $\land$ and $\lor$ (as long as both sides of the law in question 
   exist),
\item if $i$ and $j$ have no upper bound and $j'$ is arbitrary, then \underline{either} neither $i, j'$ nor 
   $i, j\land j'$ have an upper bound \underline{or} $i \lor j' = i \lor (j \land j')$, and
\item if $j,j'$ have no upper bound, then $(i \land j) \lor (i \land j') = i$ for any $i$.  \hfill $\dashv$
\end{enumerate}
\end{defin}

The following is easy to see:

\begin{lem}[{see~\cite[Observation 5]{shattered}}]
\label{almost-lattice-basic}
Assume $\la I,  \leq \ra$ satisfies (i) to (iv). Also assume $I$ has no maximal element, let $\ell \notin I$, $L = I \cup \{ \ell \}$, and
stipulate $i \leq \ell$ for all $i \in I$. Then the following are equivalent:
\begin{itemize}
\item $I$ satisfies (v) and (vi).
\item The distributive laws hold in $L$, that is, $L$ is a lattice.
\end{itemize}
\end{lem}

In fact, the only reason for having conditions (v) and (vi) in Definition~\ref{almostlattice-def} is that we want
preservation of distributivity once we add a top element to a distributive almost-lattice.

We are ready for the definition of the amalgamated limit.

\begin{defin}[{see~\cite[Definition 3]{shattered}}]   \label{amallim-def}
Given a distributive almost-lattice $\la I , \leq \ra$ and a system $( \AA_i : i \in I )$
of cBa's with complete embeddings $id : \AA_i \embed \AA_j$ for $i < j$ such that projections in
all diagrams of the form
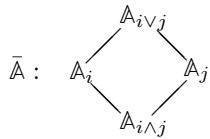
\[
\begin{picture}(50,50)(0,0)
\put(28,3){\makebox(0,0){$\AA_{i\land j}$}}
\put(-7,23){\makebox(0,0){$ \bar \AA : \;\;\; \AA_i $}}
\put(48,23){\makebox(0,0){$\AA_j$}}
\put(28,43){\makebox(0,0){$\AA_{i\lor j}$}}
\put(6,26){\line(1,1){13}}
\put(8,18){\line(1,-1){11}}
\put(31,6){\line(1,1){13}}
\put(33,38){\line(1,-1){11}}
\end{picture}
\]
are correct, we define the {\em amalgamated limit} $\AA_\ell = \AA_\amal = \lim \amal_{i\in I} \AA_i$ as follows. First
set $\AA = \bigcup_{i\in I} \AA_i \sem \{ \zero \}$.
\begin{itemize}
\item The set $D$ of nonzero conditions consists of unordered pairs $(p,q) \in \AA \times \AA$
   (so $(p,q) = (q,p)$) such that there are $i,j \in I$ with $p \in \AA_i$, $q \in \AA_j$ and
   $h_{i, i \land j} (p) = h_{j, i\land j} (q)$.
\item The order is given by stipulating $(p',q') \leq (p,q)$ if either 
\begin{itemize}
   \item $h^{i'}_{ i' \land i} (p')  \leq_{i} p$ and $h^{j'}_{ j' \land j} (q') \leq_{j} q$ or
   \item $h^{j' }_{j' \land i} (q' ) \leq_{i} p$ and $h^{i'}_{ i' \land j} (p') \leq_{j} q$ or
   \item $h^{i' }_{ i' \land i} (p')  \leq_{i} p$ and $h^{i'}_{ i' \land j} (p') \leq_{j} q$ or
   \item $h^{j' }_{ j' \land i} (q' ) \leq_{i} p$ and  $h^{j'}_{ j' \land j} (q') \leq_{j} q$
\end{itemize}
   where $(i',j')$ witnesses $(p',q') \in D$.
\end{itemize}
$\AA_\ell$ is the cBa generated by $D$, i.e., $\AA_\ell = \ro (D)$.  \hfill $\dashv$
\end{defin}

Of course, the intention is that $\AA_\ell$ is ``larger" than all the $\AA_i$ in the sense that they all completely embed into $\AA_\ell$
so that we can use this concept to define an iterated forcing construction. This is indeed the case, see Theorem~\ref{limamal-cBa} below.
Arguably, the definition of $D$ -- and thus of its completion $\AA_\ell$ -- looks complicated, in particular with respect to the order.
However, that we need to consider several cases in its definition simply has to do with the fact that our conditions are unordered pairs.
Also the use of projections in this definition is due to our considering the most general pattern for witnesses of the conditions
$(p,q)$ and $(p',q')$. In practice, when $(p',q') \leq (p,q)$ holds  we will have $i' \geq i$ and $j' \geq j$ where $(i,j)$ and $(i',j')$
witness $(p,q) \in D$ and $(p',q') \in D$, respectively, if we are in the first of the four cases defining the order above. Note that then $(p',q') \leq (p,q)$ is equivalent
to
\begin{itemize}
   \item[--] $p'  \leq_{i'} p$ and $q' \leq_{j'} q$
\end{itemize}
because $h^{i'}_{ i} (p')  \leq_{i} p$ and $p'  \leq_{i'} p$ are clearly equivalent (and similarly for $j$ and $j'$). Analogously we can simplify the
other three cases. 

One of the first things one needs to check is that the order $\leq$ is transitive on $D$; this is easy, see~\cite[Observation 6]{shattered}
for details. However, the order on $D$ is in general not separative, that is, distinct  members of $D$  may be equivalent as conditions
in the completion $\AA_\ell$ (see~\cite[Observation 9]{shattered} for several situations when this occurs). Furthermore, we identify
each $\AA_i$ with a subset of $D$ -- and thus with a subalgebra of $\AA_\ell$ -- via the map $p \mapsto (p,p)$. Note in this context that
for arbitrary conditions $(p,q) \in D$ we always have, by Definition~\ref{amallim-def}, that $(p,q) \leq (p,p) = p$ and $(p,q) \leq (q,q) = q$.

\begin{theorem}[{see~\cite[Main Lemma 7 and Lemma 8]{shattered}}]    \label{limamal-cBa}
If $\la I , \leq \ra$ and $( \AA_i : i \in I )$ are as in Definition~\ref{amallim-def}, then all $\AA_i$, $i \in I$, completely embed into the
amalgamated limit $\AA_\ell  = \lim \amal_{i\in I} \AA_i$. Furthermore, projections in all diagrams of the extended system $( \AA_i : i \in I \cup \{ \ell \} )$
are still correct.
\end{theorem}

We do not give a proof here because this is explained in detail in the reference~\cite{shattered}. However, we shall later prove a more general
theorem in a topological context, Theorem~\ref{limamal-top} in Section 3, from which Theorem~\ref{limamal-cBa} can be derived.

Condition (iii) in Definition~\ref{almostlattice-def} looks rather stringent, and one may wonder whether we can amalgamate if we relax this
condition. Counterexample 8 in~\cite{shattered} shows that this is not the case in general, though there may be situations when this can be done.
Of course, we could still define $\AA_\ell$ in  such a more general context, but complete embeddability (Theorem~\ref{limamal-cBa}) may not hold anymore,
and the definition would thus be useless for iteration theory.


\section{Examples for the amalgamated limit}

After two basic examples (Examples~\ref{exam1} and~\ref{exam2}) explaining in what sense the amalgamated limit generalizes the two-step amalgamation and the
direct limit, we will stress the similarity of three quite distinct constructions from the literature, namely, matrix iterations~\cite{BS89}, iterated ultrapowers~\cite{Sh700},
and shattered iterations~\cite{shattered}, by putting all of them into the amalgamated limit framework, in Examples~\ref{exam3} to~\ref{exam5}.

\smallskip

\begin{exam}[{Two-step amalgamation, see~\cite[Examples 2 and 6]{shattered}}]   \label{exam1}
Let $I = \{ 0 , 1 , 0 \land 1 \}$ be the three-element distributive almost-lattice. Let $\AA_i$ with $i \in I$ be cBa's such that $\AA_{0 \land 1} \embed
\AA_i$ for $i \in \{ 0,1 \}$. We consider $\amal_{\AA_{0 \land 1}} (\AA_0 , \AA_1) := \AA_\ell = \lim \amal_{i \in I } \AA_i$, the {\em amalgamation of $\AA_0$ and $\AA_1$ over $\AA_{0 \land 1}$}.
By definition of the amalgamated limit, we see that a dense subset of $\AA_\ell$ consists of pairs $(p,q)$ with $p \in \AA_0$, $q \in \AA_1$ and $h^0_{0 \land 1} (p) =
h^1_{0 \land 1} (q) \in \AA_{0 \land 1}$, and that the ordering is given by $(p',q') \leq (p,q)$ if $p' \leq p$ and $q' \leq q$. 

If $\AA_{0 \land 1} = \{ \zero, \one \}$ is the trivial forcing, $\AA_\ell$ is simply the product of the two forcing notions $\AA_0$ and $\AA_1$. Furthermore,
in the general case, we may construe $\AA_\ell$ as forcing first with $\AA_{0 \land 1}$ and then with the product of the quotients $\AA_0 / G_{0\land 1}$ and
$\AA_1 / G_{0 \land 1}$ where $G_{0 \land 1}$ is the $\AA_{0 \land 1}$-generic filter.

Amalgamations have not been used much in forcing theory. The best-known constructions involving them is Shelah's consistency proof of the projective
Baire property on the basis of the consistency of ZFC~\cite{Sh84}.
\hfill $\dashv$
\end{exam}

\smallskip

\begin{exam}[{Direct limit, see~\cite[Example 7]{shattered}}]    \label{exam2}
Let $I$ be a distributive lattice (so any two elements $i$ and $j$ of $I$ have a least upper bound $i \lor j$). Let $(\AA_i : i \in I )$ be a system of cBa's with complete
embeddings and correct projections (in the sense of Definition~\ref{amallim-def}). Then $\AA_\ell = \lim\amal_{i \in I} \AA_I$ is easily seen to be the usual
{\em direct limit} $\lim \dir_{i \in I} \AA_i$. The point is that if $(p,q)$ belongs to the dense set in the definition of the amalgamated limit, as witnessed by
$i$ and $j$ -- this means that $p \in \AA_i$, $q \in \AA_j$, and $h^i_{i \land j} (p) = h^j_{i \land j} (q) \in \AA_{i \land j}$ --, then $p \cdot q \in \AA_{i \lor j}$ is
a non-zero condition by correctness and clearly the pair $(p \cdot q, p \cdot q)$ (which is identified with $p \cdot q$ in $\AA_\ell$) is a condition stronger than
$(p,q)$ (in fact, it is even an equivalent condition in $\AA_\ell$, see~\cite[Observation 9 (iv)]{shattered}). Therefore conditions of the latter kind are
dense in $\AA_\ell$ -- and they constitute exactly the standard dense subset of the direct limit. 

Notice that in the special case $I$ has a maximum $i_0$, $\AA_\ell$ is simply the same as $\AA_{i_0}$.

The direct limit construction plays a fundamental role in iterated forcing theory and harks back  to Solovay and Tennenbaum's consistency proof of
Martin's axiom and Souslin's hypothesis~\cite{ST71}.
\hfill $\dashv$
\end{exam}

\smallskip

\begin{exam}[Matrix iterations]   \label{exam3}
Roughly speaking, a {\em matrix iteration} consists of a well-ordered system of parallel finite support iterations (fsi's); thus one may think of it as a two-dimensional
system of cBa's with complete embeddings between them. There are many different examples, but the general pattern is as follows. Let $\lambda$ and $\mu$ be ordinals.
By simultaneous recursion on $\alpha \leq \lambda$ define fsi's $\P^\gamma = (\PP^\gamma_\alpha , \dot \QQ^\gamma_\beta : \alpha \leq \lambda, \beta < \lambda)$ of ccc forcing, $\gamma \leq \mu$,
such that 
\begin{itemize}
\item[($\star$)] for every $\alpha$ and $\gamma < \delta \leq \mu$, $\PP_\alpha^\delta$ forces that $\dot \QQ^\gamma_\alpha$ is a subset of $\dot \QQ_\alpha^\delta$
and that every maximal antichain of $\dot \QQ^\gamma_\alpha$ belonging to $V^{\PP_\alpha^\gamma}$ is still a maximal antichain of $\dot \QQ_\alpha^\delta$.
\end{itemize}
Then one shows by induction on $\alpha \leq \lambda$ that for any $\gamma < \delta \leq \mu$, $\PP^\gamma_\alpha \embed \PP^\delta_\alpha$
(and, a fortiori, $\PP^\gamma_\alpha \embed \PP^\delta_\beta$ for any $\gamma < \delta$ and $\alpha \leq \beta$).
Indeed, this is trivial in the basic step $\alpha = 0$ (because all forcings $\PP^\gamma_0$ are trivial). In the successor step $\beta = \alpha + 1$,
$\PP^\gamma_\beta \embed \PP^\delta_\beta$ follows easily from the induction hypothesis $\PP^\gamma_\alpha \embed \PP^\delta_\alpha$ and $(\star)$.
If $\alpha$ is a limit ordinal, $\PP^\gamma_\alpha \embed \PP^\delta_\alpha$ follows from the induction hypothesis by Lemma~\ref{correctness-limdir}.
The point is that by the way fsi's are defined, correctness trivially holds and therefore the lemma can be applied. (Indeed, if $\alpha < \beta$ and $\gamma < \delta$,
and if $p^\gamma_\beta = ( p^\gamma_\alpha, \dot q^\gamma_{[\alpha, \beta)})  \in \PP^\gamma_\beta = \PP_\alpha^\gamma \star \dot \PP_{[\alpha,\beta)}^\gamma$
(where $\dot \PP^\gamma_{[\alpha,\beta)}$ is the $\PP^\gamma_\alpha$-name for the remainder forcing)
then $h^{\beta,\gamma}_{\alpha,\gamma} ( p^\gamma_\beta) = p^\gamma_\alpha$ and if we now have a condition $p^\delta_\alpha \in \PP^\delta_\alpha$
with $h^{\alpha,\delta}_{\alpha,\gamma} (p^\delta_\alpha) = p^\gamma_\alpha$, then $( p^\delta_\alpha, \dot q^\gamma_{[\alpha, \beta)})  \in \PP^\delta_\beta$
clearly is a common extension, as required.)

Let us now shift the perspective: suppose $\beta \leq\lambda$ is a limit ordinal and $\delta \leq \mu$ is arbitrary. Suppose all $\PP^\gamma_\alpha$
have been defined for $\gamma \leq \delta$ and $\alpha \leq \beta$ such that at least one of the inequalities is strict, and we want to define $\PP^\delta_\beta$.
Note that \[ I = \{ (\alpha,\gamma) : \gamma \leq \delta \mbox{ and }\alpha \leq \beta \mbox{ and (either }\gamma < \delta \mbox{ or }\alpha < \beta)\}\] is a distributive
almost-lattice. Also recall from the previous paragraph that all projections are correct. Thus we may set \[ \PP^\delta_\beta = \lim\amal_{(\alpha,\gamma) \in I} 
\PP_\alpha^\gamma.\] While this is formally different from the standard definition of $\PP^\delta_\beta$, namely \[ \PP^\delta_\beta = \lim \dir_{\alpha < \beta}
\PP^\delta_\alpha, \] it is easy to see that the two partial orders are forcing equivalent. Obviously, $\lim \dir_{\alpha < \beta} \PP^\delta_\alpha$ is a subset of
$\lim\amal_{(\alpha,\gamma) \in I}  \PP_\alpha^\gamma$. 
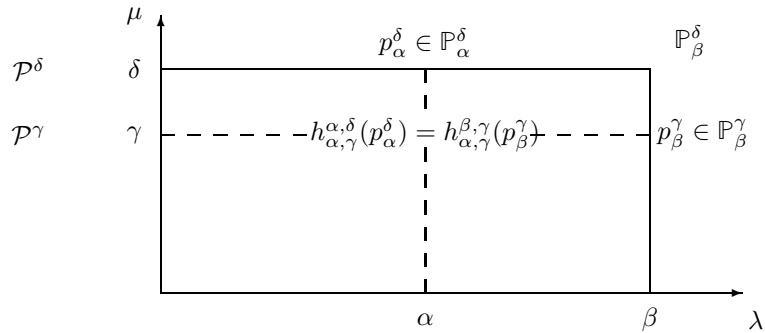
\begin{gather*}
\begin{picture}(240,140)(0,0)
\put(20,20){\vector(1,0){220}}
\put(20,20){\vector(0,1){105}}
\put(205,20){\line(0,1){85}}
\multiput(120,20)(0,10){6}{\line(0,1){5}}
\multiput(120,90)(0,10){2}{\line(0,1){5}}
\put(20,105){\line(1,0){185}}
\multiput(20,80)(10,0){6}{\line(1,0){5}}
\multiput(160,80)(10,0){5}{\line(1,0){5}}
\put(245,10){\makebox(0,0){$\lambda$}}
\put(205,10){\makebox(0,0){$\beta$}}
\put(120,10){\makebox(0,0){$\alpha$}}
\put(120,115){\makebox(0,0){$p^\delta_\alpha \in \PP^\delta_\alpha$}}
\put(10,125){\makebox(0,0){$\mu$}}
\put(10,105){\makebox(0,0){$\delta$}}
\put(10,80){\makebox(0,0){$\gamma$}}
\put(-30,105){\makebox(0,0){$\P^\delta$}}
\put(-30,80){\makebox(0,0){$\P^\gamma$}}
\put(225,80){\makebox(0,0){$p^\gamma_\beta \in \PP^\gamma_\beta$}}
\put(220,115){\makebox(0,0){$ \PP^\delta_\beta$}}
\put(120,80){\makebox(0,0){$h^{\alpha,\delta}_{\alpha,\gamma} (p^\delta_\alpha) = h^{\beta,\gamma}_{\alpha,\gamma} (p^\gamma_\beta)$}}
\end{picture} 
\end{gather*}
On the other hand, given a condition $(p^\delta_\alpha , p^\gamma_\beta)$ in the dense set defining
the amalgamated limit, that is, $p^\delta_\alpha \in \PP^\delta_\alpha$, $p^\gamma_\beta \in \PP^\gamma_\beta$, and $h^{\alpha,\delta}_{\alpha,\gamma} (p^\delta_\alpha) =
h^{\beta,\gamma}_{\alpha,\gamma} (p^\gamma_\beta) \in \PP^\gamma_\alpha$ where $\alpha < \beta$ and $\gamma < \delta$, we see that since $\PP^\gamma_\beta$ is the direct limit of the
$\PP^\gamma_{\alpha '}$, $\alpha ' < \beta$, there is $\alpha ' < \beta $ with $\alpha' \geq \alpha$ such that $p^\gamma_\beta \in \PP^\gamma_{\alpha '}$. 
This means, however, that the product $p^\delta_\alpha \cdot p^\gamma_\beta$ belongs to $\PP^\delta_{\alpha '}$, and this product clearly is a strengthening of
 $(p^\delta_\alpha , p^\gamma_\beta)$ (in fact, the two conditions are equivalent by~\cite[Observation 9 (iv)]{shattered}).

While this gives us a natural description of matrix iterations in the amalgamated limit framework, it is also clear that this approach is redundant.

Matrix iterations were originally introduced in work of Blass and Shelah~\cite{BS89} for showing the consistency of $\uu < \dd$ (in this case the $\dot \QQ^\gamma_\alpha$
are Mathias forcing with some ultrafilter $\dot \U^\gamma_\alpha$, with $\dot \U^\delta_\alpha$ carefully extending $\dot \U^\gamma_\alpha$ for $\gamma < \delta$
so that complete embeddability $(\star)$ is guaranteed) and have been used since for a plethora of consistency results; see e.g. our joint work with Fischer~\cite{BF11} which is to some extent
dual to~\cite{BS89}, the consistency of singular splitting number $\sss$ due to Dow and Shelah~\cite{DS18}, and a construction of Fischer, Friedman,
Mej\'ia, and Montoya~\cite{FFMM18}, which uses three-dimensional matrices in a special case.
\hfill $\dashv$
\end{exam}

\smallskip

\begin{exam}[{Iterating iterations, see~\cite[Section 1]{survey}}]  \label{exam4}
We now come to a construction which is in a sense orthogonal to matrix iteration -- and for which interesting examples with non-trivial amalgamated limit exist.
Recall that in the matrix approach we use the same recursion to simultaneously build many iterations which completely embed one into the other. Let us now instead
use a recursion to produce iterations, one in each step, so that the earlier iterations embed into the later ones. We call this set-up {\em iterating iterations}. Again let $\mu$ and $\lambda$ be ordinals.
By recursion on $\alpha \leq \lambda$ define iterations $\P_\alpha = (\PP^\gamma_\alpha , \dot \QQ^\delta_\alpha : \gamma \leq \mu, \delta < \mu)$ such that
\begin{itemize}
\item[($\star_1$)] for every $\gamma$ and $\alpha < \beta\leq\lambda$, $\PP^\gamma_\alpha \embed \PP^\gamma_\beta$.
\end{itemize}
Notice we do not require that the $\P_\alpha$ have direct limits and for the general construction it does not matter what $\PP^\delta_\alpha$ for limit $\delta$ is.
In the application we have in mind, $\lim \dir_{\gamma < \delta} \PP^\gamma_\alpha \embed \PP^\delta_\alpha$ and depending on $\delta$ equality may hold
(the direct limit case) or fail. We will require that
\begin{itemize}
\item[($\star_2$)] for every $\gamma$ and $\alpha < \beta \leq \lambda$, $\PP_\beta^\gamma$ forces that $\dot \QQ^\gamma_\alpha$ is a subset of $\dot \QQ_\beta^\gamma$
and that every maximal antichain of $\dot \QQ^\gamma_\alpha$ in $V^{\PP_\alpha^\gamma}$ is still a maximal antichain of $\dot \QQ_\beta^\gamma$.
\end{itemize}
Note the similarity to $(\star)$ in Example~\ref{exam3}. Clearly ($\star_2$) is sufficient to guarantee ($\star_1$) in the successor step -- that is,
if $\PP^\gamma_\alpha \embed \PP^\gamma_\beta$, then also $\PP^{\gamma+1}_\alpha \embed \PP^{\gamma+1}_\beta$ so that the issue in ($\star_1$)
is the case $\delta \leq \mu$ is a limit ordinal. If $\beta = \alpha + 1$ is a successor, the choice of $\P_\beta$ from $\P_\alpha$ will guarantee $(\star_1)$ so that
it suffices to consider limit $\beta \leq \lambda$. 

Exactly as in the previous example we now obtain a distributive almost-lattice 
\[ I = \{ (\alpha,\gamma) : \gamma \leq \delta \mbox{ and }\alpha \leq \beta \mbox{ and (either }\gamma < \delta \mbox{ or }\alpha < \beta)\}. \] 
Assume projections in all diagrams indexed by members of $I$ are correct. Then we may define the amalgamated limit
\[ \PP^\delta_\beta = \lim\amal_{(\alpha,\gamma) \in I} \PP_\alpha^\gamma\]
and, by Theorem~\ref{limamal-cBa}, we see that ($\star_1$) still holds, i.e. $\PP^\delta_\alpha \embed \PP^\delta_\beta$ for $\alpha < \beta$,
that $\PP^\gamma_\beta \embed \PP^\delta_\beta$ for $\gamma < \delta$ so that $\P_\beta$ eventually will be an iteration, and that correctness is preserved.
For a direct proof of this (and more) in this special case, see~\cite[Lemma 13]{survey} (alternatively this can be found in~\cite[Lemma 7]{Br07}, though
this paper unfortunately has a number of minor mistakes). 
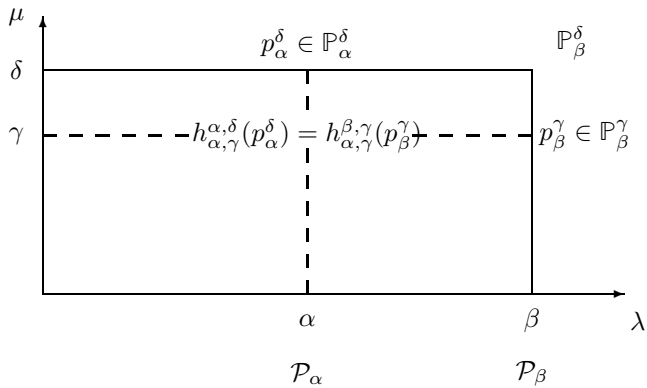
\begin{gather*}
\begin{picture}(240,160)(0,-20)
\put(20,20){\vector(1,0){220}}
\put(20,20){\vector(0,1){105}}
\put(205,20){\line(0,1){85}}
\multiput(120,20)(0,10){6}{\line(0,1){5}}
\multiput(120,90)(0,10){2}{\line(0,1){5}}
\put(20,105){\line(1,0){185}}
\multiput(20,80)(10,0){6}{\line(1,0){5}}
\multiput(160,80)(10,0){5}{\line(1,0){5}}
\put(245,10){\makebox(0,0){$\lambda$}}
\put(205,10){\makebox(0,0){$\beta$}}
\put(120,10){\makebox(0,0){$\alpha$}}
\put(205,-10){\makebox(0,0){$\P_\beta$}}
\put(120,-10){\makebox(0,0){$\P_\alpha$}}
\put(120,115){\makebox(0,0){$p^\delta_\alpha \in \PP^\delta_\alpha$}}
\put(10,125){\makebox(0,0){$\mu$}}
\put(10,105){\makebox(0,0){$\delta$}}
\put(10,80){\makebox(0,0){$\gamma$}}
\put(225,80){\makebox(0,0){$p^\gamma_\beta \in \PP^\gamma_\beta$}}
\put(220,115){\makebox(0,0){$\PP^\delta_\beta$}}
\put(120,80){\makebox(0,0){$h^{\alpha,\delta}_{\alpha,\gamma} (p^\delta_\alpha) = h^{\beta,\gamma}_{\alpha,\gamma} (p^\gamma_\beta)$}}
\end{picture} 
\end{gather*}
Of course, conditions in a dense subset of $\PP^\delta_\beta$ are again pairs $(p^\delta_\alpha , p^\gamma_\beta)$ such that
$p^\delta_\alpha \in \PP^\delta_\alpha$, $p^\gamma_\beta \in \PP^\gamma_\beta$, and $h^{\alpha,\delta}_{\alpha,\gamma} (p^\delta_\alpha) =
h^{\beta,\gamma}_{\alpha,\gamma} (p^\gamma_\beta) \in \PP^\gamma_\alpha$ where $\alpha < \beta$ and $\gamma < \delta$.

This construction has been used by Shelah~\cite{Sh700} to show the consistency of $\dd < \aa$ and of $\uu < \aa$ using the consistency of
a measurable cardinal. We briefly sketch the $\dd < \aa$ case (see ~\cite[Section 1]{survey} for details): let $\kappa$ be measurable, let $\D$ be a $\kappa$-complete
ultrafilter on $\kappa$,  and let $\lambda^\omega = \lambda > \mu > \kappa$ be regular cardinals. $\P_0$ is a finite support iteration of Hechler forcing of length $\mu$. 
For $\beta = \alpha + 1$ successor, $\P_\beta$ is the ultrapower of $\P_\alpha$ via $\D$. Then the $\dot \QQ^\gamma_\beta$, $\gamma < \mu$,
can still be construed as Hechler forcing. For limit $\beta$ do the following: if $\delta = \gamma + 1$ is successor, $\PP^\delta_\beta = \PP^\gamma_\beta
\star \dot \QQ^\gamma_\beta$ is the iteration with Hechler forcing, and if $\delta$ is limit we take the amalgamated limit described in the previous
paragraph. In the end we will have forced $\bb = \dd = \mu$ (because the Hechler generics form a scale) and $\aa = \cc = \lambda$ (because the
ultrapower destroys mad families). For the -- slightly more complicated --  consistency of $\uu < \aa$ see~\cite{Br07}.

It turns out that in this particular construction, for limit $\delta \leq \mu$ and any $\alpha$, $\PP^\delta_\alpha$ is the direct limit of the $\PP^\gamma_\alpha$ for $\gamma < \delta$
unless $cf (\delta) = \kappa$, and for limit $\beta \leq \lambda$ and any $\gamma$, $\PP^\gamma_\beta$ is the direct limit of the $\PP^\gamma_\alpha$ for
$\alpha < \beta$ unless $cf (\beta) = \omega$, so that the amalgamated limit boils down to the direct limit in many cases (see~\cite[Lemma 14]{survey} for
details). However, if $cf (\delta) = \kappa$ and $cf (\beta) = \omega$ then $\PP^\delta_\beta$ is a non-trivial amalgamated limit of the earlier forcing notions.
We also notice that in this case, the preservation of the ccc (which is trivial in direct limits) becomes a central problem (see, again, Lemma 14 of~\cite{survey}),
for even two-step amalgamations are far from preserving the ccc (see~\cite[Counterexample 10]{shattered}).
\hfill $\dashv$
\end{exam}

\smallskip

\begin{exam}[Shattered iterations~\cite{shattered}]  \label{exam5}
Roughly speaking, a {\em shattered iteration} is a matrix-style two-dimensional iteration adding random reals in limit stages instead of Cohen reals.
This means we won't have finite supports and direct limits like in Example~\ref{exam3} anymore, and the amalgamated limit becomes a necessity.
``Adding random reals in limit stages" here means that the basic forcing which embeds into the whole iteration is a large random algebra --
this is indeed dual to  the situation of finite support iterations into which a large Cohen algebra (which can be seen as the basic step of the
iteration though this is not the usual point of view) embeds. In later steps of the iteration we want to add other reals $r$ which are generic over some of the
randoms while most of the randoms will still be random over $r$ (and not just over the ground model). This is feasible by the commutativity of
random forcing.

It actually turns out that for such iterations we also need the amalgamated limit in the successor step -- though if both ordinals in our
two-dimensional system are successors this will boil down to the two-step amalgamation explained above in Example~\ref{exam1}.
It is illustrative to see why this is so. Suppose we just have two random reals, $b_0$ and $b_1$, added by the product measure algebra
$\BB_0 \star \dot \BB_1$ which is forcing equivalent (by commutativity) to $\BB_1 \star \dot \BB_0$. Now we also want to add two
Cohen reals $c_0$ and $c_1$, by $\CC_0$ and $\CC_1$, respectively, in such a way that $c_0$ is Cohen over $b_0$ and $c_1$ is
Cohen over $b_1$ but also -- and this is the point -- that $b_1$ remains random over $c_0$ and $b_0$ remains random over $c_1$.
So we need to force with $\BB_0 \star \dot \CC_0 \star \dot \BB_1$ (which is forcing equivalent to $(\BB_0 \times \CC_0) \star \dot \BB_1$)
on the one hand and with $\BB_1 \star \dot \CC_1 \star \dot \BB_0$ (equivalent to $(\BB_1 \times \CC_1 ) \star \dot \BB_0$) on the other.
We obtain the diagram
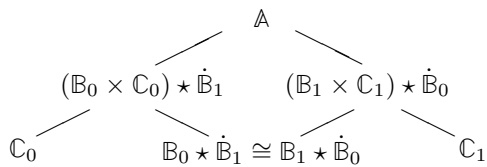
\[
\begin{picture}(180,60)(0,0)
\put(5,5){\makebox(0,0){$\CC_0$}}
\put(175,5){\makebox(0,0){$\CC_1$}}
\put(95,5){\makebox(0,0){$\BB_0 \star \dot \BB_1 \cong \BB_1 \star \dot \BB_0$}}
\put(50,30){\makebox(0,0){$(\BB_0 \times \CC_0) \star \dot \BB_1$}}
\put(135,30){\makebox(0,0){$(\BB_1 \times \CC_1) \star \dot \BB_0$}}
\put(95,55){\makebox(0,0){$\AA$}}
\put(55,37){\line(2,1){25}}
\put(55,20){\line(2,-1){20}}
\put(110,10){\line(2,1){20}}
\put(108,50){\line(2,-1){25}}
\put(10,10){\line(2,1){20}}
\put(148,20){\line(2,-1){20}}
\end{picture}
\]
in which the two-step amalgamation $\AA = \amal_{\BB_0 \star \dot \BB_1} ((\BB_0 \times \CC_0) \star
\dot \BB_1 , (\BB_1 \times \CC_1) \star \dot \BB_0 )$ of $(\BB_0 \times \CC_0) \star \dot \BB_1$ and $(\BB_1 \times \CC_1 ) \star \dot \BB_0$
over $\BB_0 \star \dot \BB_1$ arises naturally. This is the basic building block for shattered iterations.

Let us now look at the general situation, in a somewhat simplified way, in line with the two preceding examples.
Assume again $\mu$ and $\lambda$ are ordinals. Let $\mu \dot\cup \lambda$ denote the disjoint union of $\mu$ and $\lambda$.
The basic forcing $\BB_{\mu\dot\cup\lambda}$ is a measure algebra with index set $\mu \dot\cup \lambda$.
Also assume $\delta \leq \mu$ and $\beta \leq \lambda$ (this time we do not require $\delta$ and $\beta$ to be limit ordinals).
Consider again the distributive almost-lattice 
\[ I = \{ (\alpha,\gamma) : \gamma \leq \delta \mbox{ and }\alpha \leq \beta \mbox{ and (either }\gamma < \delta \mbox{ or }\alpha < \beta)\}. \] 
Assume we have constructed (by recursion on the well-founded set $I$) a system of cBa's $(\AA^\gamma_\alpha : (\alpha,\gamma) \in I)$ 
with complete embeddings $\AA^\gamma_\alpha \embed \AA^{\gamma '}_{\alpha '}$ for $\gamma \leq \gamma'$ and $\alpha \leq \alpha' $ such that 
\[ \AA^\gamma_\alpha = \EE^\gamma_\alpha \star \dot \BB_{[\gamma,\delta) \dot\cup [\alpha,\beta)} \]
where $\EE^\gamma_\alpha$ is some cBa containing $\BB_{\gamma \dot\cup \alpha}$ as a complete subalgebra
(so $\EE^\gamma_\alpha$ adds the ``first" $\gamma \dot\cup \alpha$ random reals and may also add some other generic objects,
and the ``next" $[\gamma,\delta) \dot\cup [\alpha,\beta)$ random reals will be generic over the $\EE^\gamma_\alpha$ extension), and
projections in all diagrams indexed by members of $I$ are correct. Then we let
\[ \AA^\delta_\beta = \lim \amal_{ (\alpha,\gamma) \in I} \AA^\gamma_\alpha . \]
By Theorem~\ref{limamal-cBa} we know all $\AA^\gamma_\alpha$, $(\alpha,\gamma) \in I$, completely embed into $\AA^\delta_\beta$ and
correctness is preserved.
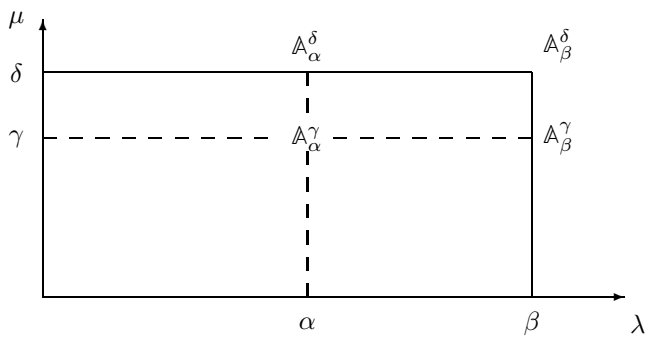
\begin{gather*}
\begin{picture}(240,140)(0,0)
\put(20,20){\vector(1,0){220}}
\put(20,20){\vector(0,1){105}}
\put(205,20){\line(0,1){85}}
\multiput(120,20)(0,10){6}{\line(0,1){5}}
\multiput(120,90)(0,10){2}{\line(0,1){5}}
\put(20,105){\line(1,0){185}}
\multiput(20,80)(10,0){9}{\line(1,0){5}}
\multiput(130,80)(10,0){8}{\line(1,0){5}}
\put(245,10){\makebox(0,0){$\lambda$}}
\put(205,10){\makebox(0,0){$\beta$}}
\put(120,10){\makebox(0,0){$\alpha$}}
\put(120,115){\makebox(0,0){$ \AA^\delta_\alpha$}}
\put(10,125){\makebox(0,0){$\mu$}}
\put(10,105){\makebox(0,0){$\delta$}}
\put(10,80){\makebox(0,0){$\gamma$}}
\put(215,80){\makebox(0,0){$\AA^\gamma_\beta$}}
\put(215,115){\makebox(0,0){$\AA^\delta_\beta$}}
\put(120,80){\makebox(0,0){$\AA^\gamma_\alpha$}}
\end{picture} 
\end{gather*}
Note we explained here only what we do in the single step $(\beta,\delta)$. To set up the whole iteration a more complicated framework,
using quadruples of ordinals instead of pairs as indices, is appropriate; see~\cite[Section 4]{shattered} for details. Correctness is not trivial here
(unlike for fsi), for we need to show that if a system has correct projections then this is still true if all cBa's of the system are iterated
with random forcing~\cite[Lemma 3]{shattered}. Last, but not least, as already hinted at in the previous example, preservation of the ccc is a central problem~\cite[Section 3]{shattered}.

Shattered iterations have been used in~\cite{shattered} to show the consistency of $\non$(meager) $> \cov$(meager) $>
\aleph_1$. This is done with the shattered iteration of Cohen forcing, that is, for each pair of ordinals $(\alpha,\gamma)$ with $\alpha < \lambda$
and $\gamma < \mu$ a Cohen real $c_{\alpha,\gamma}$ generic over $\BB_{\alpha \dot\cup\gamma}$ is added -- with the remaining
random reals still random over this Cohen real. It is fairly easy to see that this will force $\cov$(meager) $= \non$(null) $= \mu$ and 
$\cov$(null) $= \non$(meager) $= \lambda$ where $\mu < \lambda$ are regular uncountable cardinals~\cite[Section 4, Facts 23 to 26]{shattered}. 
For more results see~\cite[Section 6]{shattered}.
\hfill $\dashv$
\end{exam}


\section{The topological interpretation of the amalgamated limit}

By Stone duality, all results about (complete) Boolean algebras have a topological interpretation, and it turns out that in case of the amalgamated
limit this topological version is even more general. We first review basics of duality, and then describe the topological construction corresponding to Theorem~\ref{limamal-cBa}.

Given a Boolean algebra $\AA$ (not necessarily complete), its {\em Stone space} $X_\AA$ consists of all ultrafilters on $\AA$, with the topology
given by the elements of $\AA$: for each $a \in \AA$, $O_a = \{ x \in X_\AA : a \in x \}$ is a basic open set, and every open set is a union of such sets.
It is well-known that $X_\AA$ is a compact zero-dimensional (= there is a basis consisting of clopen sets) Hausdorff space (indeed,
since $X_\AA \sem O_a = O_{- a}$ all the $O_a$ are clopen), and, dually, given such a space $X$ there is a Boolean algebra $\AA$ such that
$X = X_\AA$, so that there is a one-to-one correspondence between Boolean algebras and compact zero-dimensional Hausdorff spaces~\cite[Chapter 3]{Ko}.
Next, if $\AA_0 \sub \AA_1$ are Boolean algebras, then the projection mapping $p = p^1_0 : X_{\AA_1} \to X_{\AA_0} , x \mapsto x \cap \AA_0$ is a continuous
surjection. Assume additionally $\AA_0$ and $\AA_1$ are cBa's. The following is well-known but we include a proof for the sake of completeness.

\begin{fact}
$\AA_0 \embed \AA_1$ iff $p$ is an open mapping (that is, images of open sets are open).
\end{fact}

\begin{proof}
Let $h= h^1_0 : \AA_1 \to \AA_0$ be the projection mapping corresponding to the embedding. For the forward direction, it clearly suffices
to show that images of basic clopen sets are clopen. Let $a_1 \in \AA_1$. We claim that $p ( O_{a_1}) = O_{h (a_1)}$. Indeed, let $x \in O_{a_1}$,
that is $a_1 \in x$. Then, by completeness of the embedding, $h(a_1) \in x$ and thus $p(x) = x \cap \AA_0 \in O_{h (a_1)}$. On the other hand, if $y \in O_{h (a_1)}$,  i.e. $h(a_1) \in y$,
then $y \cup \{ a_1 \}$ is a filter base on $\AA_1$ and therefore can be extended to an ultrafilter $x$. Clearly $p(x) = y$ and therefore $y \in p ( O_{a_1})$.

For the backward direction, assume the embedding is not complete. So there is $a_1 \in \AA_1 \sem \{ \zero \}$ such that $h(a_1) = \zero$.
Thus $\bigvee \{ a_0 \in \AA_0 : a_0 \cdot a_1 = \zero$ in $\AA_1 \} = \one$ in $\AA_0$, and $\bigcup \{ O_{a_0} : a_0 \in \AA_0 $ and $ a_0 \cdot a_1 = \zero$ in $\AA_1 \}$
is open dense in $X_{\AA_0}$. Since this set clearly is disjoint from $p (O_{a_1})$, the latter cannot be open in $X_{\AA_0}$.
\end{proof}

We next reinterpret correctness in the topological context.

\begin{fact}
Assume we have cBa's $\AA_{0 \land 1} \embed \AA_i \embed \AA_{0 \lor 1}$, $i \in \{ 0,1\}$. Projections in the diagram 
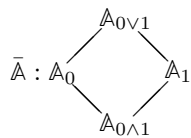
\[ 
\begin{picture}(80,50)(0,0)
\put(48,5){\makebox(0,0){$\AA_{0\land 1}$}}
\put(17,25){\makebox(0,0){$\bar \AA :  \AA_{0}$}}
\put(68,25){\makebox(0,0){$\AA_{1}$}}
\put(48,45){\makebox(0,0){$\AA_{0\lor 1}$}}
\put(26,28){\line(1,1){13}}
\put(28,20){\line(1,-1){11}}
\put(51,8){\line(1,1){13}}
\put(53,40){\line(1,-1){11}}
\end{picture}
\]
are correct iff for all $x_0 \in X_{\AA_0}$ and $x_1 \in X_{\AA_1}$ such that $p^0_{0\land 1} (x_0) = p^1_{0 \land 1} (x_1)$ there is $x_{0 \lor 1} \in X_{\AA_{0 \lor 1}}$
such that $p^{0\lor 1}_i (x_{0\lor 1} ) = x_i$ for $i \in \{ 0,1 \}$.
\end{fact}

\begin{proof}
For the forward direction, it suffices to check that $x_0 \cup x_1$ is a filter base on $\AA_{0 \lor 1}$, for we can then extend it to an ultrafilter $x_{0 \lor 1}$, which
obviously has the required properties. To see this, take $a_i \in x_i$, $i \in \{ 0,1 \}$. Then $h^i_{0 \land 1} (a_i) \in p^i_{0 \land 1} (x_i)$ (see the previous proof).
In particular $h^0_{0 \land 1} (a_0)$ and $h^1_{0 \land 1} (a_1)$ must be compatible in $\AA_{0 \land 1}$ because they belong to the same ultrafilter. 
By correctness, $a_0$ and $a_1$ are compatible in $\AA_{0 \lor 1}$, as required.

For the backward direction, let $a_i \in \AA_i$, $i \in \{ 0,1 \}$, be such that $a_{0 \land 1} := h^0_{0 \land 1} (a_0) = h^1_{0 \land 1} (a_1)$. By the argument in the previous proof,
any ultrafilter $x_{0 \land 1} \in O_{a_{0\land 1}}$ can be extended to ultrafilters $x_i \in O_{a_i}$ on $\AA_i$ and clearly $p^i_{0\land 1} (x_i) = x_{0 \land 1}$.
Therefore, by assumption, there is an ultrafilter $x_{0 \lor 1}  \in X_{\AA_{0 \lor 1}}$  with $p^{0\lor 1}_i (x_{0\lor 1} ) = x_i$. Clearly
$x_{0 \lor 1} \in O_{a_0} \cap O_{a_1}$ in $X_{\AA_{0 \lor 1}}$. This means that $a_0$ and $a_1$ must be compatible in $\AA_{0 \lor 1}$.
\end{proof}

Let us now forget about the cBa's and move to the more general context of compact Hausdorff spaces (i.e. we do not require zero-dimensionality anymore).
Then we obtain:

\begin{theorem}   \label{limamal-top}
Let $\la I , \leq \ra $ be a distributive almost-lattice, and let $( X_i : i \in I )$ be non-empty compact Hausdorff spaces. Assume that for $i \leq j$ there are continuous open
surjections $p^j_i : X_j \to X_i$ such that
\begin{itemize}
\item $i \leq j \leq k$ implies $p^k_i = p^j_i \circ p^k_j$ (commutativity) and
\item if $i \lor j$ exists then for all $x_i \in X_i$ and $x_j \in X_j$ with $p^i_{ i \land j} (x_i) = p^j_{i \land j} (x_j)$ there is $x_{i \lor j} \in X_{i \lor j}$ with
   $p^{i \lor j}_i (x_{i \lor j}) = x_i$ and $p^{i \lor j}_j (x_{i \lor j}) = x_j$ (correctness).
\end{itemize}
Then there are a compact Hausdorff space $X_\ell$ and continuous open surjections $p^\ell_i: X_\ell \to X_i$, $i \in I$, such that
\begin{itemize}
\item $i \leq j$ implies $p^\ell_i = p^j_i \circ p^\ell_j$ (commutativity) and
\item  for all $x_i \in X_i$ and $x_j \in X_j$ with $p^i_{ i \land j} (x_i) = p^j_{i \land j} (x_j)$ there is $x_{\ell} \in X_{\ell}$ with
   $p^{\ell}_i (x_{\ell}) = x_i$ and $p^{\ell}_j (x_{\ell}) = x_j$ (correctness).
\end{itemize}
\end{theorem}

\begin{proof}
Let
\[ X_\ell : = \left\{ x_\ell \in \prod_{i \in I} : p^j_i (x_\ell (j)) = x_\ell (i) \mbox{ for all } i \leq j \right\} \sub \prod_{i \in I} X_i \]
equipped with the product topology. Note that if $x \in \prod_{i \in I} X_i \sem X_\ell$ then $p^j_i (x (j)) \neq x (i)$ for some $i < j$ and since
$X_i$ is Hausdorff, there is an open neighborhood of $x$ disjoint from $X_\ell$. Therefore $X_\ell$ is a closed subset of $\prod_{i \in I} X_i $
and thus a compact Hausdorff space itself. Define the projection mapping
\[ p^\ell_i : X_\ell \to X_i, x_\ell \mapsto x_\ell (i). \]
By the definition of the product topology, $p^\ell_i$ clearly is continuous. Also, by definition of $X_\ell$ and the $p^\ell_i$, we
see
\[ p^\ell_i (x_\ell) = x_\ell (i) = p^j_i (x_\ell (j)) = p^j_i (p^\ell_j (x_\ell)) \]
so that commutativity holds. We next show that $X_\ell \neq \emptyset$ and, in fact, all $p^\ell_i$ are surjective. To this end, say that $F \sub I$ is {\em closed}
if for all $i,j \in F$, $i \land j$ and $i \lor j$ (if it exists) also belong to $F$. Clearly, if $F$ is finite, its closure $\cl (F)$ (the smallest closed set containing $F$) is finite as well.
For finite closed $F \sub I$, let 
\[ Y_F = \left\{ x \in \prod_{j \in F} X_j :  p^j_k (x(j)) = x (k) \mbox{  for all }   k \leq j \mbox{ belonging to } F \right\}.\] 

\begin{sclaim}  \label{firstclaim}
Let $F \sub G \sub I$ be finite closed sets and let $x \in Y_F$. Then there is $y \in Y_G$ extending $x$.
\end{sclaim}

\begin{proof}
Clearly it suffices to show this in case $G$ is the closure of $F \cup \{ j \}$ where $j \in I \sem F$ is arbitrary. Note that since $I$ is a distributive almost-lattice and $F$ is closed, 
either $F$ has a single maximal element $j_0$ or it has two top elements $j_0$ and $j_1$ such that $j_0 \lor j_1$ does not exist. 

In the former case, either $j_0 \lor j$ exists and we let $y (j_0 \lor j) $ be such that
$p^{j_0 \lor j}_{j_0} (y (j_0 \lor j)) = x (j_0)$ and then extend $y$ to $G$ by projecting, or $j_0 \lor j$ does not exist and there must be a second top element $j'$
of $G$ with $j' \geq j$ in which case we choose $y (j') $ such that $p^{j'}_{j_0 \land j'} ( y (j')) = p^{j_0}_{j_0 \land j'} (x (j_0))$ and we extend again by projecting.

In the latter case, either $j_0 \lor j$ or $j_1 \lor j$ exists or both exist. If only one exists, say $j_0 \lor j$, let $y ( (j_0 \lor j) \land j_1 ) = p^{j_1}_{(j_0 \lor j) \land j_1} (x (j_1))$ and note
that $p^{j_0}_{j_0 \land j_1} (x (j_0)) = p^{(j_0 \lor j) \land j_1}_{j_0 \land j_1} ( y ( (j_0 \lor j) \land j_1) )$ so that by correctness there is $y (j _0 \lor j) \in X_{j_0 \lor j}$ such that
$p^{j_0 \lor j}_{j_0} ( y (j_0 \lor j)) = x (j_0)$ and $p^{j_0 \lor j}_{  (j_0 \lor j) \land j_1 } (y (j_0 \lor j)) = y ( (j_0 \lor j) \land j_1 )$. 
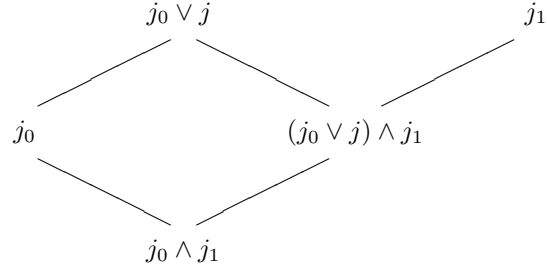
\[ 
\begin{picture}(190,100)(0,0)
\put(55,5){\makebox(0,0){$j_0 \land j_1 $}}
\put(-5,50){\makebox(0,0){$j_0$}}
\put(55,95){\makebox(0,0){$ j_0 \lor j $ }}
\put(120,50){\makebox(0,0){$ (j_0 \lor j) \land j_1 $}}
\put(190,95){\makebox(0,0){$j_1$ }}
\put(0,40){\line(2,-1){50}}
\put(60,15){\line(2,1){50}}
\put(60,85){\line(2,-1){50}}
\put(0,60){\line(2,1){50}}
\put(130,60){\line(2,1){50}}
\end{picture}
\]
Then extend $y$ to the rest of $G$ by
projecting. 

If both exist, first let $y ( j_0 \land j) = p^{j_0}_{j_0 \land j} ( x( j_0))$ and  $y ( j_1 \land j) = p^{j_1}_{j_1 \land j} ( x( j_1))$ and note that $p^{j_0 \land j}_{j_0 \land 
j_1 \land j} ( y (j_0 \land j)) = p^{j_1 \land j}_{j_0 \land j_1 \land j} ( y (j_1 \land j))$ so that by correctness we can find $y (j) \in X_j$ such that $p^j_{j_0 \land j} ( y(j)) =
y (j_0 \land j)$ and $p^j_{j_1 \land j} ( y(j)) = y (j_1 \land j)$. Next, again by correctness, find $y (j_0 \lor j) \in X_{j_0 \lor j}$ and $y (j_1 \lor j) \in X_{j_1 \lor j}$
such that $p^{j_0 \lor j}_{j_0} ( y (j_0 \lor j)) = x (j_0)$, $p^{j_0 \lor j}_{j} ( y (j_0 \lor j)) = y (j)$, $p^{j_1 \lor j}_{j_1} ( y (j_1 \lor j)) = x (j_1)$, and $p^{j_1\lor j}_{j} ( y (j_1 \lor j)) = y (j)$.
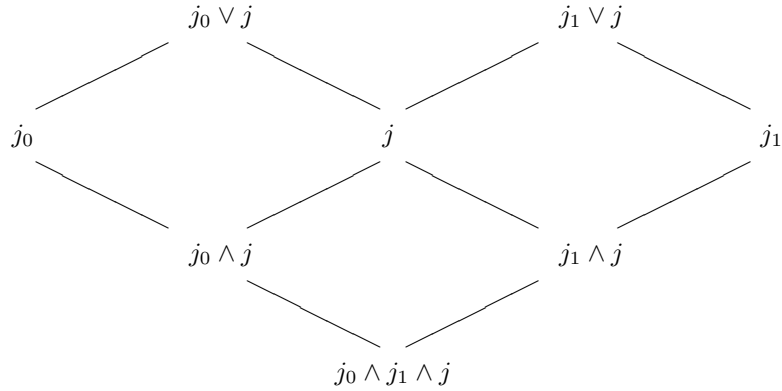
\[ 
\begin{picture}(240,150)(0,0)
\put(115,5){\makebox(0,0){$j_0 \land j_1 \land j$}}
\put(50,50){\makebox(0,0){$j_0 \land j$}}
\put(190,50){\makebox(0,0){$j_1 \land j$}}
\put(-25,95){\makebox(0,0){$j_0$}}
\put(115,95){\makebox(0,0){$ j $ }}
\put(50,140){\makebox(0,0){$ j_0 \lor j $}}
\put(260,95){\makebox(0,0){$j_1$ }}
\put(190,140){\makebox(0,0){$j_1 \lor j$}}
\put(60,40){\line(2,-1){50}}
\put(120,15){\line(2,1){50}}
\put(-20,85){\line(2,-1){50}}
\put(120,85){\line(2,-1){50}}
\put(60,60){\line(2,1){50}}
\put(200,60){\line(2,1){50}}
\put(60,130){\line(2,-1){50}}
\put(200,130){\line(2,-1){50}}
\put(-20,105){\line(2,1){50}}
\put(120,105){\line(2,1){50}}
\end{picture}
\]
Finally extend to $G$ by projecting.
\end{proof}

\begin{sclaim}  \label{secondclaim}
Let $x_i \in X_i$. There is $x_\ell \in X_\ell$ such that $x_\ell (i) = x_i$. More generally, if $F \sub I$ is closed and $y \in Y_F$, then there is $x_\ell \in X_\ell$ extending $y$.
\end{sclaim}

\begin{proof}
Clearly, the first statement is a special case of the second for $F = \{ i \}$. The general statement follows from the previous claim by compactness: fix $F$ and $y \in Y_F$ as required.
For finite closed $G \supseteq F$ define
\[ Z_G = \left\{ x \in \prod_{i \in I} X_i : y \sub x \mbox{ and } p^j_i (x(j)) = x(i) \mbox{ for all } i \leq j \mbox{ belonging to } G \right\} . \]
Clearly $Z_G$ is a closed subset of the product $\prod_{i \in I} X_i$. By Claim~\ref{firstclaim}, all $Z_G$ are non-empty. Therefore, using compactness, we see that
$\bigcap \{ Z_G : G \supseteq F$ is finite closed$\}$ is non-empty and any $x_\ell$ in this intersection belongs to $X_\ell$ and extends $y$.
\end{proof}

Thus all $p^\ell_i$ are surjective. We next show correctness. Let $i,j \in I$. If $i \lor j$ exists, this follows from correctness of the original system and the previous
claim. If $i \lor j$ does not exist, $F = \{ i, j, i \land j \}$ is a closed set, and if $x_i \in X_i$ and $x_j \in X_j$ satisfy $p^i_{i \land j} (x_i) = p^j_{i \land j} (x_j)$,
then defining $y$ by $y (i) = x_i$, $y(j) = x_j$, and $y (i \land j) = p^i_{i \land j} (x_i)$, we see that $y \in Y_F$. By the previous claim, there is $x_\ell \in X_\ell$
extending $y$, and this $x_\ell$ witnesses correctness.

We are left with showing the $p^\ell_i$ are open. To this end let $A = X_\ell \cap \left( \prod_{i \in F} U_i \times \prod_{i \in I \sem F} X_i\right) \sub X_\ell$ be a basic open set
with $F \sub I$ being finite and $U_i \sub X_i$ open for $i \in F$. 
Without loss of generality, we may assume $F$ is closed. Then $F$ has either one maximal element $i_0$ or two top elements $i_0$ and $i_1$.
In the first case let 
\[ V_{i_0} = \bigcap_{i \leq i_0} \left( p_i^{i_0} \right)^{-1} \left( U_i \right) \sub U_{i_0} \mbox{ and } V_i = p^{i_0}_i \left(V_{i_0}\right) \sub U_i \]
for $i \in F$. In the second case, first define $V_{i_n} ' = \bigcap_{i \leq i_n} \left( p_i^{i_n} \right)^{-1} \left( U_i \right) \sub U_{i_n}$ for $n \in \{ 0,1 \}$, and then
let
\[ V_{i_0 \land i_1} = p^{i_0}_{i_0 \land i_1} \left( V_{i_0} ' \right) \cap p^{i_1}_{i_0 \land i_1} \left( V_{i_1} ' \right) \sub U_{i_0 \land i_1} \]
and
\[ V_{i_n} = V_{i_n} ' \cap \left(  p^{i_n}_{i_0 \land i_1} \right)^{-1} \left( V_{i_0 \land i_1} \right) \sub U_{i_n} \mbox{ and } V_i = p^{i_n}_i \left(V_{i_n}\right) \sub U_i \]
for $i \in F$ with $i \leq i_n$, $n \in \{ 0,1 \}$ (it is easy to check that this is well-defined for $i \in F$ with $i \leq i_0 \land i_1$). We now claim:

\begin{sclaim}   \label{thirdclaim}
All $V_i$, $i \in F$, are open and for all $i \leq j$ from $F$, $V_i = p^j_i (V_j)$.
\end{sclaim}

\begin{proof}
Openness clearly follows from continuity and openness of the mappings $p^j_i$. 
Furthermore, if $j \leq i_0$, we have $V_i  = p^{i_0}_i (V_{i_0} ) = p^j_i ( p^{i_0}_j (V_{i_0} )) = p^j_i (V_j)$. Similarly if $j \leq i_1$ in the second case.
\end{proof}

\begin{sclaim}  \label{fourthclaim}
$A = X_\ell \cap \left(  \prod_{i \in F} V_i \times \prod_{i \in I \sem F} X_i\right) $.
\end{sclaim}

\begin{proof}
Since the $V_i$ are subsets of the $U_i$, the set on the right-hand side clearly is a subset of $A$. We show the converse inclusion. First consider the case
when we only have one maximal element $i_0$. Let $x_\ell \in A$. Then $p^{i_0}_i ( x_\ell (i_0)) = x_\ell (i) \in U_i$ for any $i \in F$. In particular
$x_\ell (i_0) \in ( p^{i_0}_i)^{-1} (U_i)$ for any $i \in F$, and $x_\ell ({i_0}) \in V_{i_0}$ follows. Therefore $x_\ell (i) = p^{i_0}_i ( x_\ell (i_0))  \in p^{i_0}_i (V_{i_0}) = V_i$
as well, and $x_\ell$ belongs to the right-hand side.

When we have two top elements $i_0$ and $i_1$, we basically argue in the same way. First argue that $x_\ell ( i_n) \in V_{i_n} '$ for $n \in \{ 0,1 \}$.
Thus $x_\ell (i_0 \land i_1) = p^{i_0}_{i_0 \land i_1} (x_\ell (i_0)) = p^{i_1}_{i_0 \land i_1} (x_\ell (i_1)) \in V_{i_0 \land i_1}$ and therefore
$x_\ell ( i_n) \in V_{i_n} $ for $n \in \{ 0,1 \}$. Hence $x_\ell (i) = p^{i_n}_i ( x_\ell (i_n))  \in p^{i_n}_i (V_{i_n}) = V_i$ for $i \leq i_n$, $n \in \{ 0,1 \}$,
as well, and $x_\ell$ belongs to the right-hand side.
\end{proof}

\begin{sclaim}  \label{fifthclaim}
Given $i \in F$ and $x_i \in V_i$ there is $y \in Y_F \cap \left( \prod_{j \in F} V_j \right)$ such that $y(i) = x_i$. Therefore $p^\ell_i (A) = V_i$ for all $i \in F$.
\end{sclaim}

\begin{proof}
First assume we have one maximal element $i_0$ in $F$. So $i \leq i_0$ and $V_i = p^{i_0}_i (V_{i_0})$. Thus we can find $y(i_0) \in V_{i_0}$ such
that $p^{i_0}_i (y (i_0)) = x_i$. Now extend $y$ to all of $Y_F$ by projecting. 

If we have two top elements $i_0$ and $i_1$ in $F$ we may assume without loss of generality that $i \leq i_0$. As in the previous paragraph we find
$y(i_0) \in V_{i_0}$ such  that $p^{i_0}_i (y (i_0)) = x_i$.  Let $y (i_0 \land i_1) = p^{i_0}_{i_0 \land i_1} (y (i_0)) \in p^{i_0}_{i_0 \land i_1} (V_{i_0}) =
V_{i_0 \land i_1}$. Since $ V_{i_0 \land i_1} = p^{i_1}_{i_0 \land i_1}  (V_{i_1})$, we can find $y(i_1) \in V_{i_1}$ such that 
$p^{i_1}_{i_0 \land i_1} (y (i_1)) = y (i_0 \land i_1)$. Again extend $y$ to all of $Y_F$ by projecting.

For the second statement, note that $p^\ell_i (A) \sub V_i$ is immediate by the previous claim (Claim~\ref{fourthclaim}), while $V_i \sub p^\ell_i (A)$ follows from the first part
and Claim~\ref{secondclaim}.
\end{proof}

Thus the maps $p^\ell_i$ are all open, and the proof of the theorem is complete.
\end{proof}

{\bf Acknowledgment.} We thank David Fremlin for having pointed out the topological connection to us many years ago.



\end{document}